\documentclass{amsart}

\usepackage{amsmath}
\usepackage{amsthm}
\usepackage{amsfonts}
\usepackage{amssymb}
\usepackage{psfig}
\usepackage{amscd}

\theoremstyle{plain}
\newtheorem{thm}{Theorem}
\newtheorem{lemma}[thm]{Lemma}

\newtheorem{prop}[thm]{Proposition}

\theoremstyle{definition}

\theoremstyle{remark}

\newcommand{\Q}{\ensuremath \mathbb{Q}}

\newcommand{\Z}{\ensuremath \mathbb{Z}}

\newcommand{\N}{\ensuremath \mathbb{N}}

\begin{document}

\title[OPN's Have At Least Nine Distinct Prime Factors]
{Odd Perfect Numbers Have At Least Nine Distinct Prime Factors}
\author{Pace P. Nielsen}
\address{Department of Mathematics, University of California, Berkeley, CA 94720}
\email{pace@math.berkeley.edu} 

\begin{abstract}
An odd perfect number, $N$, is shown to have at least nine distinct
prime factors. If $3\nmid N$ then $N$ must have at least twelve
distinct prime divisors.  The proof ultimately avoids previous
computational results for odd perfect numbers.
\end{abstract}

\maketitle

\section{Introduction}

A perfect number is one where $\sigma(N)=2N$.  In other words, the
sum of the divisors of $N$ is twice $N$.  These numbers have been
studied since antiquity. It is known that $N$ is an even perfect
number if and only if $N=(2^{p}-1)2^{p-1}$ with $2^{p}-1$ prime.
Sufficiency was proven by Euclid and necessity by Euler.  A prime
number of the form $2^{p}-1$ is called a Mersenne prime, and there
are currently 43 known.  There is an ongoing, online, distributed
search for such primes at \verb"http://www.mersenne.org".

The search for odd perfect numbers has not been as successful.
Currently there are none known, making their existence the oldest
unanswered question in mathematics.  However there are a great
number of necessary conditions, which go through periodic
improvements. The list of conditions given here is the same list as
in \cite{Voight}, but with recent improvements included.

Let $N$ be an odd perfect number (if such exists). Write
$N=\prod_{i=1}^{k}p_{i}^{a_{i}}$ where each $p_{i}$ is prime,
$p_{1}<p_{2}<\ldots < p_{k}$, and $k=\omega(N)$ is the number of
distinct prime factors.  The factors $p_{i}^{a_{i}}$ are called the
\emph{prime components} of $N$. Then:
\begin{itemize}
\item \emph{Eulerian Form:} We have $N=\pi^{\alpha}m^{2}$ for some
integers $\pi,\alpha, m\in\Z_{+}$, with $\pi\equiv \alpha\equiv
1\pmod{4}$ and $\pi$ prime.  The prime $\pi$ is called the
\emph{special} prime of $N$.

\item \emph{Lower Bound:} Brent, Cohen, and te Riele \cite{BCR} using a
computer search found that $N>10^{300}$.  William Lipp, using the
same techniques, is close to pushing the bound to $N>10^{500}$, and
plans to start a distributed search at the website
\verb"http://www.oddperfect.org".

\item \emph{Upper Bound:} Dickson \cite{Dickson} proved that there are finitely
many odd perfect numbers with a fixed number of distinct prime
factors. Pomerance \cite{Pomerance} gave an effective bound in terms
of $k$. This was improved in succession by Heath-Brown \cite{HB},
Cook \cite{Cook}, and finally Nielsen \cite{Nielsen} to $2^{4^{k}}$.

\item \emph{Large Factors:} Jenkins \cite{Jenkins} proved that $p_{k}>10^7$, and
Iannucci \cite{Iannucci1},\cite{Iannucci2} proved $p_{k-1}>10^4$ and
$p_{k-2}>10^2$.

\item \emph{Small Factors:}  The smallest prime factor satisfies
$p_{1}<\frac{2}{3}k+2$ as proved by Gr{\"u}n \cite{Grun}. For
$2\leqslant i\leqslant 6$, Kishore \cite{KishoreBounds} showed that
$p_{i}<2^{2^{i-1}}(k-i+1)$, and this has been slightly improved by
Cohen and Sorli \cite{CS}.

\item \emph{Number of Total Prime Factors:} Hare \cite{Hare} proved that the
total number of (not necessarily distinct) prime factors of $N$ must
be at least $47$. In unpublished work he has improved this to $75$.

\item \emph{Number of Distinct Prime Factors:} Chein \cite{Chein}
and Hagis \cite{HagisDistinct} independently proved that
$\omega(N)\geqslant 8$. Hagis \cite{Hagis3Bound} and Kishore
\cite{Kishore3Bound} showed that if $3\nmid N$ then
$\omega(N)\geqslant 11$.  This paper improves both of these bounds
by $1$.

\item \emph{The Exponents:} For the non-special primes, $p_{i}$,
write $a_{i}=2b_{i}$.  If $d=\gcd(2b_{i}+1)$ then $d \not\equiv
0\pmod{3}$ by a result of McDaniel \cite{McDaniel}.

\end{itemize}

With such a number of conditions, it might seem that an odd perfect
number could not exist.  Pomerance has given an interesting
heuristic, available at Lipp's website, suggesting that odd perfect
numbers are very unlikely.

I would like to thank William Lipp for providing some of the
factorizations given in the lemmas. Also, much of the terminology,
notation, and lemmas of the early sections of this paper match those
found in \cite{Voight} (which in turn match those in
\cite{Iannucci1}, and earlier work like \cite{Pomerance2}) in an
effort to establish a sense of both continuity and improvement, and
I'd like to thank John Voight for some interesting conversations on
topics related to his paper.

\section{Fixed Notations and Conventions}

Let $N$ be an odd perfect number and $k=\omega(N)$.  We will write
$\pi$ for the special prime.  By $\N$ we mean the non-negative
integers and by $\Z_{+}$ the positive integers. We will use an
algorithm to prove our results.  At each stage in the algorithm
there will be prime divisors of $N$ that are \emph{known} and some
that are \emph{unknown}, meaning that the prime divisors are either
specified by the algorithm or they are not, respectively. This set
of known primes will change at every stage of the algorithm as it
runs through different cases, and so the known and unknown primes
are constantly changing.  A more formal definition will be given in
a later section. We let $k_{1}$ be the number of known, distinct
prime divisors of $N$ (at any given stage), and let $k_{2}=k-k_{1}$
be the number of unknown, distinct prime divisors. Among the known
prime divisors of $N$, some of the prime \emph{components} are also
known (again, \emph{known} being a technical term meaning
\emph{specified by the algorithm}). In other words, if $p$ is a
known prime divisor of $N$ and if our algorithm yields some
$a\in\Z_{+}$ so that $p^{a}||N$, we say $p^{a}$ is a known prime
component.  We let $\ell_{1}$ be the number of known prime
components of $N$, and let $\ell_{2}=k-\ell_{1}$ be the number of
unknown prime components.

A word of warning:  In some theorems we will assume $p$ is a prime
with $p^{a}||N$, but this doesn't even mean that $p$ is a known
prime, let alone that the prime component is known. Further, even if
$p$ is a known prime, the component may still be unknown (in this
technical sense) even though we (from our position) are given that
$p^{a}||N$.  In other words, we are given the hypothesis $p^{a}||N$
but the algorithm may not have specified $p$ or $a$.  We will,
throughout, only use the words ``known'' and ``unknown'' in the
sense of known to us \emph{through our algorithm}, rather than by
hypothesis.

\section{Cyclotomic Integers}

The equation $\sigma(N)=2N$ can be (trivially) rewritten as
$\sigma(N)/N=2$ which tells us that each odd prime divisor of
$\sigma(N)$ must somehow divide $N$, and vice versa.  Thus, we want
to study the prime factorization of
\[
\sigma(N)=\prod_{i=1}^{k}\sigma(p_{i}^{a_{i}})=\prod_{i=1}^{k}\frac{p_{i}^{a_{i}+1}-1}{p_{i}-1}.
\]
Letting $\Phi_{n}(x)$ be the $n$th cyclotomic polynomial (i.e. the
minimal polynomial over $\Q$ for a primitive $n$th root of unity),
we have the partial factorization
\[
p^{n}-1=\prod_{d|n}\Phi_{d}(p)
\]
and so
\begin{equation}\label{Equation1}
\sigma(p^{n-1})=\frac{p^{n}-1}{p-1}=\prod_{d|n,\,d>1}\Phi_{d}(p).
\end{equation}
We are further interested in the factorization of $\Phi_{d}(p)$. If
$c$ and $d$ are integers with $d>1$ and $\gcd(c,d)=1$ we write
$o_{d}(c)$ for the multiplicative order of $c$ modulo $d$.  If $p$
is prime, we write $v_{p}$ for the valuation associated to $p$.  In
other words, for $n\in\Z_{+}$ we have $p^{v_{p}(n)}||n$.  The
following results of Nagell \cite{Nagell} are fundamental and are
often left as exercises in modern abstract algebra books.

\begin{lemma}[{\cite[Lemma 3]{Voight}}]\label{Lemma1}
Let $m>1$ be an integer, and let $q$ be prime.  Write $m=q^{b}n$
with $\gcd(q,n)=1$.

If $b=0$ then
\[
\Phi_{m}(x)\equiv 0\pmod{q}
\]
is solvable if and only if $q\equiv 1\pmod{m}$.  The solutions are
those $x$ with $o_{q}(x)=m$. Further,
$v_{q}(\Phi_{m}(x))=v_{q}(x^{m}-1)$ for such solutions.

If $b\neq 0$ then
\[
\Phi_{m}(x)\equiv 0\pmod{q}
\]
is solvable if and only if $q\equiv 1\pmod{n}$.  The solutions are
those $x$ with $o_{q}(x)=n$.  Further, if $m>2$ then
$v_{q}(\Phi_{m}(x))=1$  for such solutions.
\end{lemma}

We have immediately from Lemma~\ref{Lemma1} and
Equation~\ref{Equation1}:
\begin{lemma}[{\cite[Equation 4]{Iannucci1},\cite[Lemma 4]{Voight}}]\label{Lemma2}
Let $p$ and $q$ be primes, $q\geqslant 3$, and $a\in \Z_{+}$.  Then
\[
v_{q}(\sigma(p^{a}))=
\begin{cases}
v_{q}(p^{o_{q}(p)}-1)+v_{q}(a+1)& \text{ if } o_{q}(p)|(a+1) \text{
and } o_{q}(p)\neq 1,\\
v_{q}(a+1) & \text{ if } o_{q}(p)=1,\\
0 & \text{ otherwise}.
\end{cases}
\]
\end{lemma}

It turns out that the first case of Lemma~\ref{Lemma1} (when $b=0$)
is the more common means of obtaining factors for odd perfect
numbers. However, to make sure we can always reduce to that case we
need the following result, usually attributed to Bang \cite{Bang},
but given other proofs such as in \cite{BV}.

\begin{lemma}\label{Lemma3}
Let $m,x\in\Z_{+}$ with $x\geqslant 2$.  Then $\Phi_{m}(x)$ is
divisible by a prime $q$ with $o_{q}(x)=m$, except if $x=2$ and
$m=1$ or $6$, or if $x=2^{i}-1$ (for some $i\in \Z_{+}$) and $m=2$.
\end{lemma}

Note that we are only interested in this lemma when $x=p$ is a prime
dividing an odd perfect number $N$.  Thus, the case $x=2$ never
happens. Also, if $m=2$ this corresponds to the special prime (in
the Eulerian form) so $x=\pi\equiv 1\pmod{4}$ and hence cannot be of
the form $2^{i}-1$.  Thus, both exceptions in the lemma do not
affect our work.

Define $\sigma_{i}(n)=\sum_{d|n}d^{i}$, for $i\in\Z$ and
$n\in\Z_{+}$. It is clear that each of these functions is
multiplicative, $\sigma_{1}=\sigma$ is the usual sum of divisors
function, and $\sigma_{0}$ is the number of divisors function.  The
following is immediate:

\begin{lemma}\label{Lemma4}
Let $N$ be an odd perfect number.  If $p^{a}||N$, where $p$ is
prime, then for each $d|(a+1)$ the number $\Phi_{d}(p)$ is divisible
by a prime $q$ with $o_{q}(p)=d$ and $q\equiv 1\pmod{d}$.  In
particular, $\sigma(p^{a})$ has at least $\sigma_{0}(a+1)-1$
distinct prime factors in common with $N$.
\end{lemma}

\section{Fermat Primes}\label{Section:Fermat}

A prime $q$ is called a \emph{Fermat prime} if it is of the form
$q=2^{j}+1$ for some $j\in\Z_{+}$.  One can show it is necessary
that $j=2^i$ for some $i\in\N$. It is easily seen that if
$i=0,1,2,3,4$ then $2^{2^{i}}+1$ is prime (i.e.
$q=3,5,17,257,65537$), but no other Fermat primes are known. These
primes play a special role in the study of odd perfect numbers. This
is because the prime factorization of $q-1$ is exactly a power of
$2$.

The first thing we can do is restate Lemma~\ref{Lemma2} in terms of
divisors of $N$, and Fermat primes.

\begin{lemma}[{\cite{Pomerance2},\cite[Lemma 5]{Voight}}]\label{Lemma5}
Let $N$ be an odd perfect number, $p^{a}||N$ with $p$ prime, and let
$q$ be a Fermat prime.  Then:
\[
v_{q}(\sigma(p^{a}))=
\begin{cases}
v_{q}(p+1)+v_{q}(a+1) & \text{ if } \pi=p\equiv -1\pmod{q},\\
v_{q}(a+1)& \text{ if } p\equiv 1\pmod{q},\\
0 & \text{ otherwise}.
\end{cases}
\]
\end{lemma}
\begin{proof}
First suppose that $p\equiv 1\pmod{q}$.  In this case $o_{q}(p)=1$
and so we just use Lemma~\ref{Lemma2}.

Next suppose $p\not\equiv 1\pmod{q}$.  Note $o_{q}(p)$ must be a
divisor of $q-1=2^{2^{k}}$.  However, from the Eulerian form, we see
that $a+1$ is divisible by $2$ if and only if $p=\pi$.  Further, in
this case $a\equiv 1\pmod{4}$ and hence $2||(a+1)$.  Thus
$v_{q}(\sigma(p^{a}))=0$ unless $p=\pi$ and $o_{q}(p)=2$.  But this
last equality is equivalent to $p\equiv -1\pmod{q}$.
Lemma~\ref{Lemma2} gives us the needed valuation in this case.
\end{proof}

Lemma~\ref{Lemma5} tells us that there are only a few sources in
$\sigma(N)$ for copies of $q$.  In particular, if $q^{n}|N$ for some
large $n$, we can force the size of the special prime to be large.
To prove this we first need another lemma.

\begin{lemma}[{\cite[Proposition 9]{Voight}}]\label{Lemma6}
Let $N$ be an odd perfect number, let $q$ be a Fermat prime, and
suppose $p^{a}||N$ with $p$ prime.
\begin{itemize}
\item[(i)] If $p\neq\pi$ and $q^{b}|\sigma(p^{a})$ then $\sigma(p^{a})$ is
divisible by $b$ distinct primes $r_{1}, r_{2},\ldots, r_{b}$ with
$r_{i}\equiv 1 \pmod{q^{i}}$.

\item[(ii)] If $p=\pi$, $p\equiv -1\pmod{q}$, and $q^{c}|(a+1)$, then
$\sigma(p^{a})$ is divisible by $2c$ distinct primes
$r_{1},r_{1}',\ldots,r_{c},r_{c}'$ with $r_{i}\equiv r_{i}'\equiv
1\pmod{q^{i}}$.
\end{itemize}
\end{lemma}
\begin{proof}
For part (i), by Lemma~\ref{Lemma5} we have $q^{b}|(a+1)$.  So we
take $r_{i}$ to be the divisor of $\Phi_{q^{i}}(p)$ specified by
Lemma~\ref{Lemma4}. Part (ii) follows from the same lemmas, noticing
that since $p=\pi$ is special we have $2|(a+1)$, and hence we can
take $r_{i}$ and $r_{i}'$ to be the factors specified by
Lemma~\ref{Lemma4} of $\Phi_{q^{i}}(p)$ and $\Phi_{2q^{i}}(p)$
respectively.
\end{proof}

\begin{prop}[{c.f. \cite[Proposition 10]{Voight}}]\label{Prop7}
Let $N$ be an odd perfect number, and let $q$ be a Fermat prime with
$q^{n}|N$.  Suppose $k$, $k_{1}$, $k_{2}$, $\ell_{1}$, and
$\ell_{2}$ have their usual meanings. Further suppose
$q^{b}||\sigma(\text{known prime components of }N)$. Finally, let
$k_{1}'$ (respectively $\ell_{1}'$) be the number of distinct prime
factors among the $k_{1}$ known prime divisors (respectively, among
the $\ell_{1}$ known prime components) which are congruent to
$1\pmod{q}$.

If
\[
\tau=n-b-(k_{1}'-\ell_{1}'+k_{2})(k_{1}'+k_{2}-1)>0
\]
then $\pi\equiv -1\pmod{q}$ and $\pi$ is among the unknown prime
components.  If, further, each known prime, $p$, with unknown
component and with $p\equiv 1\pmod{4}$ satisfies $p\not\equiv
-1\pmod{q^{\tau}}$, then $\pi$ is among the unknown primes, and
\[
v_{q}(\pi+1)\geqslant
n-b-(k_{1}'-\ell_{1}'+k_{2}-1)(k_{1}'+k_{2}-2)-\lfloor
(k_{1}'+k_{2}-1)/2 \rfloor=\tau'.
\]
\end{prop}
\begin{proof}
First, we will show the contrapositive of the initial statement.
Suppose $\pi$ is among the known components, or if it isn't then
$\pi\not\equiv -1\pmod{q}$. Let $p^{a}$ be an unknown component of
$N$.  Lemma~\ref{Lemma5}, combined with what we've just said,
implies $v_{q}(\sigma(p^{a}))=v_{q}(a+1)$ if $p\equiv 1\pmod{q}$,
and equals $0$ otherwise. Using the equation $\sigma(N)/N=2$, there
are $n-b$ copies of $q$ that must be accounted for by $\sigma$ of
the unknown components. Note that there are at most $k_{1}'+k_{2}$
distinct prime factors of $N$ that are congruent to $1\pmod{q}$.
Thus, by Lemma~\ref{Lemma6}, at most $k_{1}'+k_{2}-1$ copies of $q$
divide $a+1$ if $p\equiv 1\pmod{q}$ (since $p$ itself cannot divide
$\sigma(p^{a})$), otherwise we end up with too many distinct prime
divisors of $N$ which are congruent to $1\pmod{q}$. There are at
most $k_{1}'-\ell_{1}'+k_{2}$ primes $p$ with unknown component
satisfying $p\equiv 1\pmod{q}$. Hence at most
$(k_{1}'-\ell_{1}'+k_{2})(k_{1}'+k_{2}-1)$ copies of $q$ can be
accounted for by $\sigma$ of the unknown components. Therefore
\[\tau=n-b-(k_{1}'-\ell_{1}'+k_{2})(k_{1}'+k_{2}-1)\leqslant 0\]
or we will have left over copies of $q$ not accounted for in
$\sigma(N)$.

To get the last statement, now suppose $\tau>0$, $\pi\equiv
-1\pmod{q}$, and $\pi$ is among the unknown components.
Lemma~\ref{Lemma5} tells us that the only place the extra $\tau$
copies of $q$ can be accounted for is $v_{q}(\pi+1)$, and hence
$q^{\tau}|(\pi+1)$. Therefore $\pi \equiv -1 \pmod{q^{\tau}}$ and so
$\pi$ is an unknown prime.

This means that when we counted the maximum number of possible
(unknown) prime divisors $\equiv 1\pmod{q}$ we included one too
many. Thus there are at most $k_{1}'+k_{2}-1$ primes dividing $N$
that are congruent to $1\pmod{q}$, and hence at most
$(k_{1}'-\ell_{1}'+k_{2}-1)(k_{1}'+k_{2}-2)$ copies of $q$ can be
accounted for by $\sigma$ of the unknown, non-special components.
For the special component $\pi^{\alpha}$, by Lemma~\ref{Lemma6} part
(ii) we see that at most $(k_{1}'+k_{2}-1)/2$ copies of $q$ can
divide $\alpha+1$, else we obtain too many factors of $N$ congruent
to $1\pmod{q}$. Since, in fact, an integer number of copies of $q$
divides $\alpha+1$, we can take the floor.  Putting this all
together, $\tau'$ copies of $q$ must divide $\pi+1$.
\end{proof}

Thus, if a large power of a Fermat prime divides $N$, we see that
$\pi+1$, and hence $\pi$, must be large.  This isn't quite good
enough to simplify our search to manageable cases.  We need a way to
find another large prime divisor of $N$.  The trick is to consider
divisors of $\sigma(q^{n})$. Suppose $q$ is a Fermat prime, with
$q^{n}||N$, and $n$ very large. Then using the previous lemma, we
can force $\pi$ to be very large. One might wonder if
$\pi|\sigma(q^{n})$. The following result speaks to this issue.

\begin{lemma}[{\cite[Lemma 1]{BCR}}]\label{Lemma8}
If $p$ and $q$ are odd primes with $p|\sigma(q^{k})$ and
$q^{m}|(p+1)$ then $k\geqslant 3m$.
\end{lemma}

So, using the terminology of Proposition~\ref{Prop7}, if $\tau>0$,
$q^{\tau}\nmid (p+1)$ for the known primes, and also $3\tau'>n$,
then $\pi\nmid\sigma(q^{n})$ by Lemma~\ref{Lemma8}. But since $n$ is
big we would expect a large prime divisor of $N$ which divides
$\sigma(q^{n})$. The following work clarifies how large a divisor we
can find for $\sigma(q^{n})$.

\begin{lemma}\label{Lemma9}
Let $p$ be an odd prime and let $q=3$ or $5$.  If $q^{p-1}\equiv
1\pmod{p^{2}}$ then either $(q,p)=(3,11)$ or $q^{o_{p}(q)}-1$ has a
prime divisor greater than $10^{13}$.  If $q=17$ then either
$(q,p)=(17,3)$ or $q^{o_{p}(q)}-1$ has a prime divisor greater than
$10^{11}$.
\end{lemma}
\begin{proof}
The papers \cite{Montgomery} and \cite{KR} give a list of $(q,p)$
for which $q^{p-1}\equiv 1\pmod{p^{2}}$ and $p<10^{13}$ with $q=3$
or $5$ (or $p<10^{11}$ with $q=17$). In the cases $(q,p)\neq
(3,11),(17,3)$ the following table gives the requisite factor of
$q^{o_{p}(q)}-1$.
\[
\begin{array}{cc|ccc|cc}
q & & & p & & & \text{Large factor of } q^{o_{p}(q)}-1\\
\hline
3 & & & 1006003 & & & 154680726732318637\\
5 & & & 20771 & & & 625552508473588471\\
5 & & & 40487 & & & 625552508473588471\\
5 & & & 53471161 & & & 50493456782731\\
5 & & & 164533507 & & & 52082118058261\\
5 & & & 6692367337 & & & 8930008316757509\\
5 & & & 188748146801 & & & 40093613041379\\
17 & & & 46021 & & & 1365581260423071390161\\
17 & & & 48947 & & & 63895279579889
\end{array}
\]
\end{proof}

For use shortly, we make the following definition. Letting $p,q$ be
odd primes, $p\neq q$, we set
\[
o'_{q}(p)=
\begin{cases}
 & \text{ if } 2\nmid o_{q}(p), \\ &  \text { or  } 4\nmid o_{q}(p)
\text{ and either }\\
o_{q}(p) &  \quad\text{ (i) }  p=\pi \text{ or }\\
& \quad \text{ (ii) } p\equiv 1\pmod{4} \text{ and } \pi \text{ is not }\\
& \quad\quad\quad\!\! \text{ among the known components } \\
0 & \text{ otherwise}.
\end{cases}
\]
In other words, $o'_{q}(p)$ is the usual order function, unless it
is impossible for both $p^{a}||N$ and $o_{q}(p)\mid (a+1)$ to hold,
due to consideration of the Eulerian form.

\begin{prop}\label{Prop10}
Let $N$ be an odd perfect number with $k$, $k_{1}$, and $k_{2}$
having their usual meanings. Suppose $q=3$ or $5$ is a known prime
divisor of $N$, $q^{n}||N$, $q\neq \pi$ and $\pi\nmid\sigma(q^{n})$.
Suppose $p_{1},\ldots, p_{k_{1}-1}$ are the other known factors of
$N$, besides $q$.  For each $i\in[1,k_{1}-1]$ define
\[
\epsilon_{i}=
\begin{cases}
0 & \text{if } o'_{p_{i}}(q)=0 \\
\max(s+t-1,1) & \text{if } o'_{p_{i}}(q)\neq 0,\,
s=v_{p_{i}}(\sigma(q^{o_{p_{i}}(q)-1})),\,\\ & \quad\text{and
$t\in\Z_{+}$ minimal so that } p_{i}^{t}>100.
\end{cases}
\]
Set $V=\prod_{i=1}^{k_{1}-1}p_{i}^{\epsilon_{i}}$.  Suppose $\pi$ is
among the $k_{2}$ unknown prime factors, and $k_{2}>1$. Finally,
assume that all unknown prime factors are greater than $100$.

If
\[
\min\left(10^{13},\left(\frac{\sigma(q^{n})}{V}\right)^{\frac{1}{k_{2}-1}},
\left(\frac{\sigma(q^{100})}{V}\right)^{\frac{1}{k_{2}-1}}\right)>1,
\]
then $\sigma(q^{n})$ has a prime divisor among the unknown primes at
least as big as the above minimum. If $q=17$ one can replace
$10^{13}$ with $10^{11}$, and then the result still holds.
\end{prop}
\begin{proof}
We do the case when $q=3$, since the other cases are similar. First
suppose $\sigma(q^{n})=(q^{n+1}-1)/(q-1)$ is at most divisible by
$p_{i}^{\epsilon_{i}}$ for the known primes, and square-free for the
unknown primes. Then since $\pi,q\nmid \sigma(q^{n})$, the largest
unknown divisor of $\sigma(q^{n})$ is at least
\[
\left(\frac{\sigma(q^{n})}{V}\right)^{\frac{1}{k_{2}-1}},
\]
unless this quantity is $\leqslant 1$ (in which case there might be
no unknown factors).

So we may assume there is some prime $p|N$, $o'_{p}(q)\neq 0$, so
that $\sigma(q^{n})$ is divisible by $p^{2}$ if $p$ is unknown, or
$p^{\epsilon+1}$ if $p$ is known (and $\epsilon$ is the
corresponding $\epsilon_{i}$), with $p$ \emph{maximal} among such
primes. By Lemma~\ref{Lemma9} we may also assume that if $p^{2}\nmid
(q^{p-1}-1)$ and $p$ is unknown then $p|(n+1)$.  (This is where
$10^{13}$ comes into the minimum.)

Thus in either case, $p^{t}|(n+1)$, where $p^{t}>100$ (taking $t=1$
if $p$ is an unknown prime). Then we have
\[
\sigma(q^{p^{t}-1})|\sigma(q^{n}).
\]
Thus it suffices to find a large divisor of $(q^{p^{t}}-1)/(q-1)$.
By Lemma~\ref{Lemma1}, $(q^{p^{t}}-1)/(q-1)$ is only divisible by
primes larger than $p$, or $p$ itself to the first power.  (In this
case, $q$ being Fermat means the quantity isn't divisible by $p$,
and we could replace $\max(s+t-1,1)$ by $s+t-1$ in the definition of
$\epsilon_{i}$.  But to keep similar notations later when we take
$q$ to be an arbitrary prime, we don't use this fact.) But then, by
the maximality condition on $p$, $(q^{p^{t}}-1)/(q-1)$ is not
divisible by more than $p_{i}^{\epsilon_{i}}$ for known primes and
the first power for all the unknown primes. So the analysis we used
in the first paragraph goes through by only changing $n+1$ to
$p^{t}$. Finally, note that $p^{t}>100$, so we have the appropriate
bound.
\end{proof}

The most useful case when we will use Proposition~\ref{Prop10} is
when $n$ is very large and $k_{1}$ is close to $k$. So, in practice,
we will usually end up with $10^{13}$ as the lower bound on a
divisor of $\sigma(q^{n})$. In fact, a lot of work could be saved if
the bounds given in \cite{KR} were improved, or there was some means
to work around square divisors.

\section{Non-Fermat Primes}\label{Section:NonFermat}

Sometimes our odd perfect number will not be divisible by a large
power of a Fermat prime, but rather by a large power of some
arbitrary prime, $q$.  Unfortunately, we don't have
Lemma~\ref{Lemma5}, and so we can't put all of the ``extra'' factors
of $q$ into $\pi+1$. This causes two problems.  First, we have to
spread the extra factors of $q$ among the unknown primes, thus
reducing the number of extra factors we have at hand, exponentially.
Second, if $p\not\equiv 1\pmod{q}$ is one of our unknown primes, we
cannot reduce to the case $q|\Phi_{2}(p)$ but rather $q|\Phi_{d}(p)$
for some (arbitrary) $d>1$, $d|(q-1)$. Thus, we need a way of
bounding the size of $p$ for which $q^{n}|\Phi_{d}(p)$. This second
problem is easily dealt with.

\begin{lemma}[{\cite[Theorem 2]{KR}}]\label{Lemma11}
Let $q$ be an odd prime.  Let $a_{1}$ be a primitive root modulo
$q$. For $r\geqslant 2$, define $a_{r}=a_{1}^{q^{r-1}}\pmod{q^{r}}$.
Then $\{a_{r}^{m}\pmod{q^{r}}\,|\, m=0,\ldots, q-2\}$ gives a
complete set of incongruent solutions to $a^{q-1}\equiv
1\pmod{q^{r}}$.
\end{lemma}

\begin{lemma}\label{Lemma12}
Let $q<1000$ be an odd prime.  If $p^{q-1}\equiv 1\pmod{q^{n}}$ for
some $n\in\Z_{+}$ and some odd prime $p$, then $p\geqslant
\min(q^{n-2},10^{50})$ except when
\[
(p,q)=(40663372766570611389846294355914421,7).
\]
\end{lemma}
\begin{proof}
Let $q<1000$ be an odd prime.  Fix $m\in\N$ so that
$q^{m-2}>10^{50}$ and $m$ is minimal.  A computer search using
Lemma~\ref{Lemma11} demonstrates that every positive solution to
$x^{q-1}\equiv 1\pmod{q^{m}}$ with $x\neq 1$ satisfies $x>10^{50}$.
Hence the same is true if we replace $m$ with a larger integer.
Thus, it suffices to search for prime solutions to $x^{q-1}\equiv
1\pmod{q^{n}}$ for $n\leqslant m$.  Again, a computer search yields
the stated result.
\end{proof}

The choice of $1000$ in Lemma~\ref{Lemma12} is easily improved, but
is large enough for our needs.  Also, the exceptional $(p,q)$ in the
lemma is irrelevant to our work, since in this case we find
$o_{q}(p)=6$ and $\sigma(p^{5})$ gives rise to at least $12$
additional prime factors of $N$ besides $p$ and $q$.

\begin{lemma}\label{Lemma13}
Let $N$ be an odd perfect number, let $q$ be an odd prime, and
suppose $p^{a}||N$ with $p$ prime.
\begin{itemize}
\item[(i)] If $p\equiv 1\pmod{q}$ and $q^{b}|\sigma(p^{a})$ then $\sigma(p^{a})$ is
divisible by $b$ distinct primes $r_{1}, r_{2},\ldots, r_{b}$ with
$r_{i}\equiv 1 \pmod{q^{i}}$.

\item[(ii)] If $p\not\equiv 1\pmod{q}$, $o_{q}(p)|(a+1)$, and $q^{c}|(a+1)$, then
$\sigma(p^{a})$ is divisible by $c\cdot \sigma_{0}(o_{q}(p))$
distinct primes divisors $r_{i,j}$, $i\in[1,c]$,
$j\in[1,\sigma_{0}(o_{q}(p))]$, with $r_{i,j}\equiv 1\pmod{q^{i}}$.
\end{itemize}
\end{lemma}
\begin{proof}
Analogous to Lemma~\ref{Lemma6}.
\end{proof}

\begin{prop}\label{Prop14}
Let $N$ be an odd perfect number, and let $q<1000$ be a prime
divisor of $N$ with $q^{n}|N$.  Suppose $b$, $k$, $k_{1}$, $k_{2}$,
$\ell_{1}$, $\ell_{2}$, $k_{1}'$, and $\ell_{1}'$, have the same
meanings as in Proposition~\ref{Prop7}.  Suppose further that the
exceptional case of Lemma~\ref{Lemma12} doesn't hold.  Let $T$ be
the set of known primes with unknown component, different from $q$,
and $\not\equiv 1\pmod{q}$. Let
\begin{align*}
\tau=n-b & -\sum_{p\in T,\, o'_{q}(p)\neq
0}\left(v_{q}(p^{o_{q}(p)}-1)+ \left\lfloor
\frac{k_{1}'+k_{2}}{\sigma_{0}(o_{q}(p))} \right\rfloor\right) \\ &
- (k_{1}'-\ell_{1}'+k_{2})(k_{1}'+k_{2}-1).
\end{align*}

If $\tau>0$ then one of the unknown primes is not congruent to
$1\pmod{q}$.  Further, in this case, one of the unknown primes is at
least as large as $\min(q^{\tau'-2},10^{50})$ where
\begin{align*}
\tau'=  \min_{1\leqslant m\leqslant k_{2}}\Bigg{\lceil}\Bigg{(} n-b-
\sum_{p\in T,\, o'_{q}(p)\neq 0}\left(v_{q}(p^{o_{q}(p)}-1)+
\left\lfloor \frac{k_{1}'+k_{2}-m}{\sigma_{0}(o_{q}(p))}
\right\rfloor\right)\\
-(k_{1}'-\ell_{1}'+k_{2}-m)(k_{1}'+k_{2}-m-1)-m\left\lfloor
\frac{k_{1}'+k_{2}-m}{2} \right\rfloor \Bigg{)}/m\Bigg{\rceil}.
\end{align*}
\end{prop}
\begin{proof}
The proof is similar to Proposition~\ref{Prop7}.  From the equation
$\sigma(N)/N=2$ we know that $q^{n}|\sigma(N)$, and so we try to
account for as many copies of $q$ in $\sigma(N)$ as we can. The
quantity $\tau$ is exactly how many copies of $q$ are unaccounted
for, if all of the unknown primes are $\equiv 1\pmod{q}$, and we try
to account for as many copies of $q$ as possible from the known
primes using Lemmas~\ref{Lemma2} and~\ref{Lemma13}. If $\tau>0$ this
means we actually have left over copies of $q$, which yields a
contradiction, and hence not all of the unknown primes are $\equiv
1\pmod{q}$.

In this case, let $m$ be the number of unknown primes $\not\equiv
1\pmod{q}$. Lemma~\ref{Lemma2} tells us that there are two sources
of copies of $q$ in $\sigma(N)$; the exponents of the primes, and
$p^{o_{q}(p)}-1$.  The quantity
\begin{multline}\label{Eq:Long}
n -b -\sum_{p\in T,\, o'_{q}(p)\neq 0}\left(v_{q}(p^{o_{q}(p)}-1)+
\left\lfloor \frac{k_{1}'+k_{2}-m}{\sigma_{0}(o_{q}(p))}
\right\rfloor\right)\\
-(k_{1}'-\ell_{1}'+k_{2}-m)(k_{1}'+k_{2}-m-1)-m\left\lfloor
\frac{k_{1}'+k_{2}-m}{2} \right\rfloor
\end{multline}
is the number of copies of $q$ in $\sigma(N)$ not yet accounted for,
after we account for as many copies of $q$ as we can (again using
Lemmas~\ref{Lemma2} and~\ref{Lemma13}) except those copies of $q$
which come from $\sigma(p^{o_{q}(p)-1})$ for the unknown primes $p$
with $p\not\equiv 1\pmod{q}$.  Thus, if we divide
Equation~\ref{Eq:Long} by $m$, and take the ceiling, we have a
number of copies of $q$ that must be accounted for by
$p^{o_{q}(p)}-1$ for an unknown prime $p$ with $o_{q}(p)\neq 1$.
Since $m$ is unspecified by the hypotheses, we take the minimum over
all possibilities. Finally, we apply Lemma~\ref{Lemma12}.
\end{proof}

Next we want to prove a result analogous to
Proposition~\ref{Prop10}, except for arbitrary primes. However,
there have to be a few differences. First, we need something to play
the role of Lemma~\ref{Lemma8}. It turns out that it doesn't hurt
much to just assume that the large prime, call it $p$, coming from
Proposition~\ref{Prop14} may in fact divide $\sigma(q^{n})$. So we
can write $q^{n+1}-1=(q-1)p^{c}m$ for some $c\in \N$ and
$m\in\Z_{+}$. Powering this equation to the $(q-1)$st power, since
$p^{q-1}\equiv 1\pmod{q^{\tau'}}$ and $\tau'\leqslant n$, we have
\[
((q-1)m)^{q-1}\equiv ((q-1)mp^{c})^{q-1} = (q^{n+1}-1)^{q-1}\equiv
1\pmod{q^{\tau'}}
\]
and hence we can bound $(q-1)m$ using an analogue of
Lemma~\ref{Lemma12} (where we look for solutions to the equation
$x^{q-1}\equiv 1\pmod{q^{n}}$ which are divisible by $q-1$, rather
than prime solutions). The following lemma provides this result:

\begin{lemma}\label{Lemma15}
Let $q<1000$ be an odd prime.  Suppose $a^{q-1}\equiv 1\pmod{q^{n}}$
for some $n\in\Z_{+}$ and some integer $a$ with $(q-1)|a$.  If
$q\geqslant 11$ then $a\geqslant \min(q^{n-2},10^{50})$.  If $q=7$,
then $a\geqslant \min(q^{n-3},10^{50})$.
\end{lemma}
\begin{proof}
An easy computer search as in Lemma~\ref{Lemma12}.
\end{proof}

Next we need to work with possible square factors of
$\sigma(q^{n})$, similar to Lemma~\ref{Lemma9}.

\begin{lemma}\label{Lemma16}
Let $p$ and $q$ be primes with $p\in(10^2,10^{11})$ and $q=7$, $11$
or $13$. If $q^{p-1}\equiv 1\pmod{p^{2}}$ then
$\sigma(q^{o_{p}(q)-1})$ is divisible by two primes greater than
$10^{11}$.
\end{lemma}
\begin{proof}
By \cite{KR}, there are only $3$ pairs $(p,q)$ satisfying the
conditions, namely $$(p,q)=(491531,7),(863,13),(1747591,13).$$  In
the first case, since $65|o_{p}(q)$, just factor $\sigma(7^{64})$ to
find two distinct primes greater than $10^{11}$.  In the second case
$o_{p}(q)=862$.  One factor, $16002623839393$, is easily found and
another is
$$\frac{13^{431}+1}{14\cdot 863^{3}\cdot 68099}.$$
In the last case, since $195|o_{p}(q)$, one factors
$\sigma(13^{38})$ to find one prime greater than $10^{11}$ and
factors $\sigma(13^{64})$ to find the other one. (Or, to save time,
use the online factorizations on Richard Brent's extensions of the
Cunningham tables for these last two.)
\end{proof}

We are now ready to prove:

\begin{prop}\label{Prop17}
Let $N$ be an odd perfect number and let $q=11$ or $13$ be a known
prime divisor of $N$, with $q^{n}||N$. Let $\tau, \tau'$ be as in
Proposition~\ref{Prop14}, suppose all the hypotheses of that
proposition are met, and let $p$ be the guaranteed unknown prime.
Let $p_{1},\ldots,p_{k_{1}-1}$ be the known primes different from
$q$. Let $\epsilon_{i}$ be defined as before, and put
$V=\prod_{i=1}^{k_{1}-1}p_{i}^{\epsilon_{i}}$. Suppose $k_{2}>1$.
Finally, assume that all unknown prime factors are greater than
$100$.

If
\[
\min\left(10^{11},\left(\frac{10^{50}}{(q-1)V}\right)^{\frac{1}{k_{2}-1}},
\left(\frac{q^{\tau'-2}}{(q-1)V}\right)^{\frac{1}{k_{2}-1}}
\right)>1
\]
then $\sigma(q^{n})$ has a prime divisor, different from $p$, among
the unknown primes, at least as big as the above minimum. If $q=7$
the same result holds if we replace $q^{\tau'-2}$ by $q^{\tau'-3}$.
\end{prop}
\begin{proof}
We do the case $q=11$ or $13$, the other being similar. So assume
the above minimum is $>1$. First note that if $d|(n+1)$ then
$\sigma(q^{d-1})|\sigma(q^{n})$ and so it suffices to show that
$\sigma(q^{d-1})$ has a prime divisor larger than the above minimum,
different from $p$, for some $d|(n+1)$.

By the method of proof given in Proposition~\ref{Prop10}, we may
assume that at most $\epsilon_{i}$ copies of $p_{i}$ divide
$\sigma(q^{d-1})$, for some $d$ either greater than $100$ or equal
to $n+1$. Further, because $10^{11}$ occurs in the above minimum, we
may assume that the only unknown prime that may divide
$\sigma(q^{n})$ which may be greater than $10^{11}$ is $p$. Then by
Lemma~\ref{Lemma16} and the fact that the unknown primes are greater
than $100$, we may assume $\sigma(q^{n})$ is square-free for unknown
primes, except possibly $p$.

From the proof for Proposition~\ref{Prop14}, we have $p^{q-1}\equiv
1\pmod{q^{\tau'}}$.  Write $q^{d}-1=(q-1)mp^{c}$ with $c\in\N$,
$m\in\Z_{+}$, and $\gcd(p,m)=1$.  Powering this equation to the
$(q-1)$st power, we have
\[
((q-1)m)^{q-1}\equiv ((q-1)mp^{c})^{q-1} = (q^{d}-1)^{q-1}\equiv
1\pmod{q^{\min(\tau',100)}}.
\]
By Lemma~\ref{Lemma15}, $(q-1)m\geqslant \min(q^{\tau'-2}, q^{98},
10^{50})=\min(q^{\tau'-2},10^{50})$. Thus $$\frac{m}{V}\geqslant
\min\left(\frac{q^{\tau'-2}}{(q-1)V},\frac{10^{50}}{(q-1)V}\right).$$
Since $m/V$ is at least as big as the part of $\sigma(q^{d-1})$ made
up from the unknown primes, different from $p$, if we take the
$k_{2}-1$ root of the minimum we have the appropriate lower bound.
\end{proof}

\section{Abundance and Deficiency}\label{Section:AD}

Let $n\in \Z_{+}$.  Recall the multiplicative function
$\sigma_{-1}(n)=\sum_{d|n}d^{-1}$ we introduced earlier.  This
function can alternatively be written using the formula
$\sigma_{-1}(n)=\sigma(n)/n$, and so $\sigma(n)/n=2$ if and only if
$\sigma_{-1}(n)=2$. A number $n$ is called \emph{abundant} when
$\sigma_{-1}(n)>2$ and \emph{deficient} when $\sigma_{-1}(n)<2$.  We
can use abundance and deficiency computations to limit choices on
possible prime factors of an odd perfect number $N$.  First, we
extend the definition of $\sigma_{-1}$ by setting
\[
\sigma_{-1}(p^{\infty})=\lim_{a\rightarrow
\infty}\sigma_{-1}(p^{a})=\frac{p}{p-1}.
\]
\begin{lemma}[{\cite[Proposition 2]{Voight}}]\label{Lemma18}
Let $p$ and $q$ be odd primes.  If $1\leqslant a<b\leqslant \infty$
then $1<\sigma_{-1}(p^{a})<\sigma_{-1}(p^{b})$.  If $1\leqslant
a,b\leqslant \infty$ and $p<q$ then
$\sigma_{-1}(q^{b})<\sigma_{-1}(p^{a})$.
\end{lemma}

\begin{lemma}\label{Lemma19}
Let $N$ be an odd perfect number.  Suppose $p_{1},\ldots, p_{k_{1}}$
are the known prime factors of $N$, $p_{i}^{a_{i}}|N$, and
$k_{1}<k=\omega(N)$. If
$\Pi=\prod_{i=1}^{k_{1}}\sigma_{-1}(p_{i}^{a_{i}})< 2$, then the
smallest unknown prime is $$p_{k_{1}+1}\geqslant
\frac{\Pi}{2-\Pi}.$$
\end{lemma}
\begin{proof}
We find
\[
2=\sigma_{-1}(N)\geqslant
\left(\prod_{i=1}^{k_{1}}\sigma_{-1}(p_{i}^{a_{i}})\right)\sigma_{-1}(p_{k_{1}+1})
=\Pi\cdot\frac{p_{k_{1}+1}+1}{p_{k_{1}+1}}.
\]
where the inequality in the middle follows from Lemma~\ref{Lemma18}.
Noting $\Pi\geqslant 1$, we obtain $ \frac{2}{\Pi}\geqslant
1+\frac{1}{p_{k_{1}+1}}$. Therefore $\frac{2-\Pi}{\Pi}\geqslant
\frac{1}{p_{k_{1}+1}}$ and taking reciprocals gives us the result,
since $2-\Pi>0$.
\end{proof}

Note that in the lemma if $\Pi>2$ then
$\prod_{i=1}^{k_{1}}p_{i}^{a_{i}}$ is abundant, hence $N$ is
abundant. If $\Pi=2$, then $\prod_{i=1}^{k_{1}}p_{i}^{a_{i}}$ is
already an odd perfect number.

The following lemma is the true key to our search for odd perfect
numbers, as simple as the proof is (after wading through the
hypotheses). This is because we built up machinery in the last few
sections to find bounds for large prime divisors of $N$.

\begin{lemma}[{c.f. \cite[Lemma 2.2]{CS}}]\label{Lemma20}
Let $N$ be an odd perfect number.  Let $p_{1},\ldots p_{k}$ be the
prime divisors of $N$, and let $a_{i}$ be such that
$p_{i}^{a_{i}}||N$.  Fix the numbering on the indices so that
$p_{1},\ldots, p_{\ell_{1}}$ are the primes with known prime
component, $p_{\ell+1},\ldots p_{k_{1}}$ are the other known primes,
and $p_{k_{1}+1}<\ldots < p_{k}$ are the unknown primes. Suppose
among the unknown primes we have bounds $p_{k}>P_{1},\ldots,
p_{k-v+1}>P_{v}$, with $v<k_{2}$ and $P_{i}>1$ for each $i\in[1,v]$.
For each $u\in[0,v]$, set
\[
\Delta_{u} =
\left(\prod_{i=1}^{\ell_{1}}\sigma_{-1}(p_{i}^{a_{i}})\right)
\left(\prod_{i=\ell_{1}+1}^{k_{1}}\frac{p_{i}}{p_{i}-1}\right)
\left(\prod_{i=1}^{u}\frac{P_{i}}{P_{i}-1} \right).
\]
Finally, suppose $k_{2}>0$.

If $\Delta_{u}<2$ then the smallest unknown prime is
$$p_{k_{1}+1}\leqslant \frac{\Delta_{u}(k_{2}-u)}{2-\Delta_{u}}+1.$$  Therefore,
$$p_{k_{1}+1}\leqslant \min_{u\in[0,v],\,\Delta_{u}<2}\left(
\frac{\Delta_{u} (k_{2}-u)}{2-\Delta_{u}}+1\right).$$
\end{lemma}
\begin{proof}
We compute
\begin{align*}
2& =\sigma_{-1}(N)=\prod_{i=1}^{k}\sigma_{-1}(p_{i}^{a_{i}})\\
& \leqslant
\left(\prod_{i=1}^{\ell_{1}}\sigma_{-1}(p_{i}^{a_{i}})\right)
\left(\prod_{i=\ell_{1}+1}^{k-u}\sigma_{-1}(p_{i}^{\infty})\right)
\left(\prod_{i=k-u+1}^{k}\sigma_{-1}(P_{k-i+1}^{\infty}) \right)
\\ & = \Delta_{u}
\prod_{i=k_{1}+1}^{k-u}\sigma_{-1}(p_{i}^{\infty}) \leqslant
\Delta_{u}\prod_{i=0}^{k-u-k_{1}-1}\sigma_{-1}((p_{k_{1}+1}+i)^{\infty})\\
& =\Delta_{u}
\frac{p_{k_{1}+1}+k-u-k_{1}-1}{p_{k_{1}+1}-1}=\Delta_{u}\left(1+\frac{k_{2}-u}{p_{k_{1}+1}-1}\right).
\end{align*}
Now, recall that $u\leqslant v<k_{2}$ which implies
$k_{2}-u\geqslant 1$. Also $0<\Delta_{u}<2$, so we solve the main
inequality as we did in the previous lemma, finding
\[
p_{k_{1}+1}\leqslant \frac{\Delta_{u}(k_{2}-u)}{2-\Delta_{u}}+1.
\]
The last statement follows.
\end{proof}

One major difference between Lemma~\ref{Lemma20} and
Lemma~\ref{Lemma19} is that if $\Delta_{u}>2$ then that doesn't
necessarily imply $N$ is abundant.  (It is true that if $k_{1}=k$
and $\Delta_{0}<2$ then $N$ is deficient, however.)  This means that
we might end up with $\Delta_{0}>2$, and hence we have no upper
bound on $p_{k_{1}+1}$.

\section{The Algorithm}

To verify that an odd perfect number has at least 9 distinct prime
factors, we use a factor chain algorithm to check all possible cases
for odd perfect numbers with exactly 8 distinct prime factors, and
find contradictions in each case, similar to the algorithms
described in both \cite{CS} and \cite{Voight}. This section will
describe the algorithm, and the next section will explain my
implementation.

We start by knowing $3|N$, so we consider each of the cases
$3^{2}||N$, $3^{4}||N$, and so forth. We think of each of these
cases as branches on a tree.  On the branch $3^{2}||N$ we know that
$\sigma(3^{2})|2N$ and so $13|N$. Hence the case $3^{2}||N$ branches
further and we have to consider each of the subcases $13||N$,
$13^{2}||N$, $13^{4}||N$, and so forth.  We continue finding new
factors, and the tree continues branching out.

Some cases do not result in new primes.  To continue the factor
chain in these cases we use the bounds of Section~\ref{Section:AD}
to find an interval for the smallest unknown prime.  Then we
consider all the new primes in the given interval.  We only stop
branching when we reach a contradiction (such as having too many
factors).

As it stands there are still an infinite number of branches. To get
around this problem, we can combine all the cases $p^{n}||N$ for
large $n$ together into one case.  More precisely, let $B$ be a
large integer (around the size $10^{50}$) which we fix at the
beginning of the algorithm.  Then, once we reach the branch $p^{n}$,
with $n$ large enough so that $p^{n}>B$, we no longer consider the
case $p^{n}||N$ but rather we just assume $p^{a}|N$ for some
$a\geqslant n$, and thus consider all the remaining cases together.
(In particular, $p$ is still a known prime divisor of $N$, but the
component is unknown, although we are given a lower bound of $B$ on
the size of the component.)  We then label this composite branch by
$p^{\infty}$.

So, for example, if $B=10^2$ then since $3^{4}<B<3^{6}$ the first
level of our tree would consist of the branches $3^{2}$, $3^{4}$,
and $3^{\infty}$. When we are on a branch with $p^{\infty}$ we say
$p$ is an {\it infinite prime} (not to be confused with the infinite
primes of algebraic number theory).  Notice that infinite primes do
not provide more primes for the factor chain, and so we have to rely
on on the interval bounds of Section~\ref{Section:AD}. If we set $B$
too low then the primes become infinite too quickly and we may have
the case the there is no upper bound for an interval. (This
corresponds to the case when $\Delta_{0}>2$ in Lemma~\ref{Lemma20}.)
If we make $B$ large enough this never happens (for a proof see
\cite{CS}), and the algorithm will only have to consider a finite
number of cases.

We are now able to expressly define what we mean by \emph{known}
primes and components.  Suppose we are on the branch
$3^{\infty}5^{4}$.  Then the known primes are $3,5,11,71$ (the
primes $11$ and $71$ come from $\sigma(5^{4})|2N$), and the only
known component is $5^{4}$.  In this case we say that $3$ and $5$
are \emph{on} while $11$ and $71$ are \emph{off}.  In other words,
the \emph{on} primes are exactly the known primes for which we have
started the branching process.  Note that $k_{1}-\ell_{1}$ is
exactly the number of (known) primes which are infinite or off.

In our example, since $11$ is the smallest off prime we continue the
branching process on this prime, rather than $71$.  When there are
no off primes we use the interval bounds to arrive at another prime,
as explained earlier. Whenever we reach a contradiction, we go to
the next available branch.

The following is a possible, initial print-out in the case $k=5$,
$B=50$.  (Note: When calculating intervals this output only used
$\Delta_{0}$ in Lemma~\ref{Lemma20}.)
\newline
\newline
\noindent\verb"3^2 => 13^1"
\newline\noindent\verb"   13^1 => 2^1 7^1"
\newline\noindent\verb"      7^2 => 3^1 19^1"
\newline\noindent\verb"         19^"$\infty$\verb" : 21 < p_5 < 23 N"
\newline\noindent\verb"      7^"$\infty$\verb" : 9 < p_4 < 21"
\newline\noindent\verb"         11^"$\infty$\verb" : 374 < p_5 < 540"
\newline
\newline
The letter \verb"N" means that there are no primes in the interval,
which is a contradiction.  Adding more contradictions, or using the
full power of Lemma~\ref{Lemma20}, can further simplify the output.

\section{An Implementation}

There are three main differences between our implementation of this
algorithm, and the implementations in \cite{CS} and \cite{Voight}:

First, the bound $B$ is not allowed to increase within the
algorithm.  Allowing $B$ to vary fully automates the algorithm at
the expense of unnecessary complexity.  The number $B$ is fixed at
the outset, and only increased manually if needed.

Second, the use of Lemma~\ref{Lemma20} allows for stronger upper
bounds on intervals. In the terminology of that lemma, we always
have $v\leqslant 3$.  We take $P_{1}=10^{7}$, $P_{2}=10^{4}$, and
$P_{3}=10^{2}$, unless we already have known prime divisors larger
than these bounds, or unless these bounds are superseded by the work
in Sections~\ref{Section:Fermat} and~\ref{Section:NonFermat}. (So,
for example, if the largest known prime is $>10^{7}$, and the next
largest known prime is $<10^{2}$, then we could take $P_{1}=10^{4}$
and $P_{2}=10^{2}$.  If large powers of small primes divide $N$ then
we can use the work in previous sections to increase $P_{1}$ and
$P_{2}$.)

Third, the contradictions used are different.  Here is a complete
list of the contradictions in our implementation:
\begin{itemize}
\item[MT] There are too many total factors.

\item[MS] There are too many copies of a single prime with known component.

\item[S] There is an off prime smaller than an on prime coming from
interval computations.

\item[A] The number is abundant.

\item[D] There are $k$ known primes, and $\Delta_{0}<2$, hence $N$ is deficient.

\item[F] The special prime $\pi$ belongs to a known component, but the hypotheses
of Proposition~\ref{Prop7} hold showing $\pi$ must be in an unknown
component due to a Fermat prime.

\item[N] There are no primes in the interval given by Lemmas~\ref{Lemma19} and~\ref{Lemma20},
or all the known primes are on and the only primes in the interval
are already known.

\item[SF1] There are $k-1$ known primes, and the interval formula gives
an upper bound of $p_{k}<C$, but we know from the fact that a large
power of a small Fermat prime divides $N$ that some unknown prime is
larger than $C$, by Proposition~\ref{Prop7}.

\item[SF2] Similar to SF1 except we have a
contradiction between the interval formula and
Proposition~\ref{Prop10}.

\item[SNF1] Similar to SF1, except we have a contradiction from a
small non-Fermat prime, using Proposition~\ref{Prop14}.

\item[SNF2] Similar, using Proposition~\ref{Prop17}.

\item[P1Int] There are $k-1$ known primes all smaller than
$10^{7}$, and the interval formula gives an upper bound
$p_{k}<10^7$.

\item[P2Int] There are $k-2$ known primes all smaller than $10^4$,
or $k-1$ known primes with the largest $>10^7$ and none other
$>10^4$; and the interval formula says the next prime is $<10^4$.

\item[P3Int] There are $k-3$ known primes all smaller than $10^2$,
or $k-2$ known primes only one greater than $10^4$ and none other
$>10^2$, or $k-1$ known primes with one larger than $10^6$ and
another larger that $10^4$ and none other $>10^2$; and the interval
formula says the next prime is $<10^2$.
\end{itemize}
Of course, before the algorithm checks for any contradictions it
always checks if we have an odd perfect number.  There were other
contradictions we might have included.  For example, if there are
$k$ known primes, and none of them can be the special prime, this
means $N$ cannot be an odd perfect number.  However, the strength of
our method lies in the fact that we rarely have exactly $k$ known
primes, and in those cases one of the other contradictions will do.

\section{Points of Improvement and the Results}

The algorithm was implemented in Mathematica on a Pentium 4 personal
computer. The factorizations of $\sigma(p^{a})$ were carried out
using a only a probable primality test, good for integers
$<10^{16}$. Thus, if a factor of $\sigma(p^{a})$ was larger than
this bound, Mathematica's primality proving routine was also run.
After running the algorithm it became clear that certain
modifications would speed up the process. First, the bound $B$ was
too uniform.  So, $B$ was replaced by two bounds $B_{1}$ and
$B_{2}$, where a prime $p<1000$ became infinite when $p^{a}|N$ and
$p^{a}>B_{1}$, while a prime $p>1000$ became infinite when $p^{a}|N$
and $p^{a}>B_{2}$.

With $B_{1}=10^{50}$ and $B_{2}=10^{18}$, the algorithm didn't
terminate until we reached the case
$$3^{\infty}5^{\infty}17^{\infty}257^{\infty}65537^{\infty}$$
because $B_{2}$ wasn't large enough.  After increasing it to
$B_{2}=10^{30}$ the algorithm finished the case $\omega(N)=8$
without finding an odd perfect number, running for a total of less
than 3 days, after numerous stops and starts.

One of the longest cases (besides the one mentioned above) was
$3^{4}7^{\infty}5^{1}$.  This is because the bound given in
Proposition~\ref{Prop17} depends on taking $(k_{2}-1)$st roots.  So,
when $k_{2}=3$ the bound is really around the square-root of where
it should be.  This led to formulating the improved result:

\begin{prop}\label{Prop21}
Assume all of the hypotheses of Proposition~\ref{Prop10} hold.
Further, suppose that we can bound $p_{k_{1}+1}<P$.  Finally assume
$k_{2}>2$. If
\[
\min\left(10^{13},\left(\frac{\sigma(q^{n})}{V
P}\right)^{\frac{1}{k_{2}-2}}, \left(\frac{\sigma(q^{100})}{V
P}\right)^{\frac{1}{k_{2}-2}}\right)>1,
\]
(replacing $10^{13}$ by $10^{11}$ if $q=17$) then $\sigma(q^{n})$
has a prime divisor among the unknown primes at least as big as the
above minimum.
\end{prop}
\begin{proof}
One reduces to the case that $\sigma(q^{n})$ is square-free for
unknown primes, just as in Proposition~\ref{Prop10}.  Then the
result is obvious.
\end{proof}

Of course, Proposition~\ref{Prop17} is similarly improved.  Once
these improvements were put into place, the case
$3^{4}7^{\infty}5^{1}$ took less than an hour.  The case
$\omega(N)=8$ took less than half a day of non-stop running, with
less than $800,000$ lines of output. We thus have:

\begin{thm}\label{Theorem22}
An odd perfect number has at least $9$ distinct prime divisors.
\end{thm}

Running the algorithm with $B_{1}=10^{50}$, $B_{2}=10^{30}$, $k=11$,
and forcing $3\nmid N$, also terminated without finding any odd
perfect numbers after about one million lines of output.  Thus:

\begin{thm}\label{Theorem23}
An odd perfect number $N$ with $3\nmid N$ has at least $12$ distinct
prime divisors.
\end{thm}

At this point, the results seemed to depend on a lot of previous
computation.  To avoid this, the algorithm was rerun with the
following changes:
\begin{itemize}
\item All uses of the P1Int, P2Int, and P3Int contradictions
were removed.

\item The results of \cite{Iannucci1},
\cite{Iannucci2}, and \cite{Jenkins} were not used in interval
calculations.

\item It wasn't assumed that $3|N$.  In fact,
Gr{\"u}n's result wasn't even used, but rather Lemma~\ref{Lemma20},
to obtain $p_{1}\leqslant k+1$.

\item It wasn't assumed that $\omega(N)\geqslant 8$.  So
each $k\in[1,8]$ was checked.
\end{itemize}
With these changes, the only results (outside this paper) which were
used were (i) some of the easy to prove non-computational lemmas
(such as the Eulerian form) and (ii) the main result of \cite{KR}
for primes $a\leqslant 17$.  The algorithm again terminated after
about one million lines of output, verifying
Theorem~\ref{Theorem22}.  A similar run verified
Theorem~\ref{Theorem23}.

\section{Future Results}

To do the case $\omega(N)=9$ with these techniques, in a reasonable
amount of time, it is necessary to find bounds for $p_{k-2}$, or
increase the lower bound on $p_{k-1}$. For example, consider the
case of
$$3^{\infty}5^{\infty}17^{\infty}257^{\infty}65537^{\infty}.$$  We can
make the powers on these primes so large that we may, for all
intents and purposes, assume $p_{k}$ is as big as we like (except
not quite big enough to use \cite{Nielsen}). However, the best we
can do with our tools for $p_{k-1}$ is $10^{13}$, which is much
smaller than the lower bound already given by the interval formula.
The interval for $p_{k-1}$ has a length of about $10^{18}$, which is
much too large to check one prime at a time.  Further, it is
currently computationally unfeasible to push the bound in \cite{KR}
up to $10^{19}$.

Another possible line of attack would be to suppose that two primes
$q_{1}, q_{2}$ divide $N$, each to a high power.  Then one might be
able to prove that $\sigma(q_{1}^{n_{1}})$ and
$\sigma(q_{2}^{n_{2}})$ must each have large divisors, different
from one another, once $n_{1}$ and $n_{2}$ are large enough. Another
alternative would be to show that the square-free part of
$\sigma(q^{n})$ becomes large as $n$ does.

\bibliographystyle{amsplain}
\bibliography{oddperfect_bib}

\end{document}